# Lyapunov Conditions for Uniform Asymptotic Output Stability and a Relaxation of Barbălat's Lemma

**Iasson Karafyllis[*] and Antoine Chaillet[**]**

[*]Dept. of Mathematics, National Technical University of Athens,
Zografou Campus, 15780, Athens, Greece.
emails: iasonkar@central.ntua.gr , iasonkaraf@gmail.com

[**]L2S-CentraleSupélec and Université Paris-Saclay,
Institut Universitaire de France.
email: antoine.chaillet@centralesupelec.fr

**Abstract**
Asymptotic output stability (AOS) is an interesting property when addressing control applications in which not all state variables are requested to converge to the origin. AOS is often established by invoking classical tools such as Barbashin-Krasovskii-LaSalle's invariance principle or Barbălat's lemma. Nevertheless, none of these tools allow to predict whether the output convergence is uniform on bounded sets of initial conditions, which may lead to practical issues related to convergence speed and robustness. The contribution of this paper is twofold. First, we provide a testable sufficient condition under which this uniform convergence holds. Second, we provide an extension of Barbălat's lemma, which relaxes the uniform continuity requirement. Both these results are first stated in a finite-dimensional context and then extended to infinite-dimensional systems. We provide academic examples to illustrate the usefulness of these results and show that they can be invoked to establish uniform AOS for systems under adaptive control.

## 1. Introduction

Several control applications require to impose a prescribed behavior not necessarily to the whole state, but only to specific state variables or some output of interest. A typical example is that of adaptive control, in which it is not necessary to impose a zero steady-state error on the parameter estimation to get a proper behavior of the closed-loop system. Many notions are available for the description of such a property, including partial stability [22], which specifically addresses the case when only part of the variables are requested to have a stable behavior, and output stability [6,17], which considers stability of a specific targeted output.

  Beyond convergence of the output to the origin, a desirable additional property is that the rate at which this convergence is achieved is uniform on bounded sets of initial states. This leads to the notion of Uniform Asymptotic Output Stability (UAOS). This uniformity in initial states comes for free in the full-state stability analysis of finite-dimensional systems (in the sense that stability and attractivity imply uniform convergence rate, [8]), but is not guaranteed in the context of output stability [14].

UAOS can be casted in the more general framework of stability with respect to two measures [10,13], meaning a $KL$ estimate in which the left-hand side measure is the output norm and the right-hand side measure is the full state norm. This expression captures the fact that transient behavior of the output usually depends on the full initial state, and not only on the initial value of the considered output. UAOS is then fully characterized in terms of Lyapunov functions [6,11,21], which allows in turn to analyze robustness to exogenous inputs [18].

From an application perspective, uniformity is a desirable property as it precludes the possibility to have arbitrarily slow convergence of the output from a bounded set of initial states. But, perhaps more importantly, lack of uniformity may have detrimental consequences on the robustness to exogenous inputs or model uncertainties. For instance, an example was given in [14] of a (non-uniformly) asymptotically output-stable system whose output does not converge to zero in the presence of an arbitrarily small vanishing input (although the system evolves on a bounded set). Due to the above-mentioned Lyapunov characterizations, such lack of robustness is precluded by UAOS. Indeed, Theorem 1 in [21] indicates that a $KL$ stability estimate (i.e., UAOS) holds if and only if there exists a Lyapunov function.

Classical tools to establish that a specific output converges to zero are the well-known Barbashin-Krasovskii-LaSalle's invariance principle and Barbălat's lemma [8]. Both these results help guaranteeing that the output $y(t)$ converges to zero under the assumption that a Lyapunov-like function $V$ satisfies, along the system's solutions,

$$\dot{V}(x) \leq -\rho(W(x)), \text{ with } W(x) \geq a(|y|) \tag{1.1}$$

where $\rho$ is a continuous positive definite function and $a$ is a function of class $K_\infty$. Nevertheless, none of these two results guarantees that the output convergence is indeed uniform. An example was actually given in [14] of a system satisfying (1.1) but whose solutions do not uniformly converge to the origin.

Our first main result in this paper is to show that condition (1.1) does ensure UAOS provided that $W$ does not increase along the system's solution (Theorem 1). We show how this result can be useful to establish semiglobal uniform asymptotic output stability through adaptive control for systems satisfying the matching condition [9]. More precisely we show that, given any compact set of initial conditions, there exists an adaptive control law that makes the closed-loop system UAOS from this set (Theorem 3).

We also provide novel conditions to ensure (possibly non-uniform) asymptotic output stability by relaxing the uniform continuity requirement imposed by Barbălat's lemma to merely uniform upper semi-continuity from the right and uniform lower semi-continuity from the left (Theorem 2). We also provide a testable criterion to ensure this "quasi-uniform continuity" property in practice. While other extensions of Barbălat's lemma were already given the literature [4,19,20], we are not aware of any existing results allowing to relax the uniform continuity requirement. A consequence of the two main contributions can be roughly summarized as follows: the differential inequality (1.1) ensures asymptotic output stability if $W$ does not increase or decrease more than linearly in time along the system's solutions and the convergence rate is uniform if $W$ does not increase along solutions.

We then extend all these results in an infinite-dimensional setup that allows the study of partial differential equations and time-delay systems. In the particular case of time-delays systems, our results generalize Lyapunov-Krasovskii conditions for UAOS presented in [5] as they do not require a dissipation rate involving the whole Lyapunov-Krasovskii functional itself, which often proves hard to show in practice.

**Notation.** Throughout this paper, we adopt the following notation.

* $\Re_+ := [0, +\infty)$.

* Let $S \subseteq \Re^n$ be an open set and let $A \subseteq \Re^n$ be a set that satisfies $S \subseteq A \subseteq cl(S)$, where $cl(S)$ is the closure of $S \subseteq \Re^n$. By $C^0(A; \Omega)$, we denote the class of continuous functions on $A$, which take values in $\Omega \subseteq \Re^m$. By $C^k(A; \Omega)$, where $k \geq 1$ is an integer, we denote the class of functions on $A \subseteq \Re^n$, which takes values in $\Omega \subseteq \Re^m$ and has continuous derivatives of order $k$. In other words, the functions of class $C^k(A; \Omega)$ are the functions which have continuous derivatives of order $k$ in $S = \text{int}(A)$ that can be continued continuously to all points in $\partial S \cap A$. When $\Omega = \Re$ then we write $C^0(A)$ or $C^k(A)$.

* By $K$ we denote the class of strictly increasing $C^0$ functions $a: \Re_+ \to \Re_+$ with $a(0) = 0$. By $K_\infty$ we denote the class of strictly increasing $C^0$ functions $a: \Re_+ \to \Re_+$ with $a(0) = 0$ and $\lim_{s \to +\infty} a(s) = +\infty$. By $KL$ we denote the set of all continuous functions $\sigma: \Re_+ \times \Re_+ \to \Re_+$ with the properties: (i) for each $t \geq 0$ the mapping $\sigma(\cdot, t)$ is of class $K$; (ii) for each $s \geq 0$, the mapping $\sigma(s, \cdot)$ is non-increasing with $\lim_{s \to +\infty} \sigma(s, t) = 0$.

* For a vector $x \in \Re^n$, $|x|$ denotes its Euclidean norm and $x'$ denotes its transpose.

* Let $S \subseteq \Re^n$ be a non-empty set with $0 \in S$. We say that a function $V: S \to \Re_+$ is positive definite if $V(x) > 0$ for all $x \in S$ with $x \neq 0$ and $V(0) = 0$. We say that a continuous function $V: S \to \Re_+$ is radially unbounded if the following property holds: "for every $M > 0$ the set $\{x \in S : V(x) \leq M\}$ is compact". For $V \in C^1(S; \Re_+)$ we define $\nabla V(x) = \left(\frac{\partial V}{\partial x_1}(x), \ldots, \frac{\partial V}{\partial x_n}(x)\right)$.

## 2. Output Stability for Finite-Dimensional Systems

<u>2.1. Definitions</u>

We start by introducing the class of dynamical system considered here and by recalling notions of output stability. Let $S \subseteq \Re^n$ be a non-empty set with $0 \in S$, $f: S \to \Re^n$ be a locally Lipschitz vector field with $f(0) = 0$ and $h: S \to \Re^k$ be a continuous mapping with $h(0) = 0$. Consider the dynamical system

$$\dot{x} = f(x) \qquad \qquad (2.1)$$
$$x \in S$$

with output

$$y = h(x) \qquad \qquad (2.2)$$

We assume that the dynamical system (2.1) is forward complete, i.e., for every $x_0 \in S$ the unique solution $x(t) = \phi(t, x_0)$ of the initial-value problem (2.1) with initial condition $x(0) = x_0 \in S$ exists for all $t \geq 0$ and satisfies $\phi(t, x_0) \in S$ for all $t \geq 0$.

We use the notation $y(t, x_0) = h(\phi(t, x_0))$ for all $t \geq 0$, $x_0 \in S$ and $B_R = \{ x \in \mathfrak{R}^n : |x| < R \}$ for all $R > 0$. The following properties are standard in the analysis of output stability: see for instance [6,18,21,22].

**Definition 1 (Output Stability Notions):** *Let $\Omega \subseteq S$ be a positively invariant set for (2.1) with $0 \in \Omega$. We say that system (2.1), (2.2) is*
  i) *Lagrange output stable on $\Omega$ if for every $R > 0$ the set $\{ |y(t, x_0)| : x_0 \in B_R \cap \Omega, t \geq 0 \}$ is bounded.*
  ii) *Lyapunov output stable on $\Omega$ if for every $\varepsilon > 0$ there exists $\delta(\varepsilon) > 0$ such that for all $x_0 \in B_{\delta(\varepsilon)} \cap \Omega$, it holds that $|y(t, x_0)| \leq \varepsilon$ for all $t \geq 0$.*
  iii) *Asymptotically Output Stable (AOS) on $\Omega$ if system (2.1), (2.2) is Lagrange and Lyapunov output stable on $\Omega$ and $\lim_{t \to +\infty} (y(t, x_0)) = 0$ for all $x_0 \in \Omega$.*
  iv) *Uniformly Asymptotically Output Stable (UAOS) on $\Omega$ if system (2.1), (2.2) is Lagrange and Lyapunov output stable on $\Omega$ and for every $\varepsilon, R > 0$ there exists $T(\varepsilon, R) > 0$ such that for all $x_0 \in B_R \cap \Omega$, it holds that $|y(t, x_0)| \leq \varepsilon$ for all $t \geq T(\varepsilon, R)$.*

The following two lemmas provide characterizations of the stability notions of Definition 1 in terms of comparison functions. Since Lemmas 1 and 2 are special cases of Lemmas 4 and 5 (given below) their proofs are omitted.

**Lemma 1 (Output stability through a $K_\infty$ estimate):** *Let $\Omega \subseteq S$ be a positively invariant set for (2.1) with $0 \in \Omega$. System (2.1), (2.2) is Lagrange output stable on $\Omega$ and Lyapunov output stable on $\Omega$ if and only if there exists a function $\zeta \in K_\infty$ such that the following estimate holds for all $x_0 \in \Omega$ and $t \geq 0$:*

$$|y(t, x_0)| \leq \zeta(|x_0|)$$

**Lemma 2 (UAOS through a $KL$ estimate):** *Let $\Omega \subseteq S$ be a positively invariant set for (2.1) with $0 \in \Omega$. System (2.1), (2.2) is UAOS on $\Omega$ if and only if there exists a function $\sigma \in KL$ such that the following estimate holds for all $x_0 \in \Omega$ and $t \geq 0$:*

$$|y(t, x_0)| \leq \sigma(|x_0|, t)$$

2.2. Problem Statement

Powerful mathematical results to establish AOS are Barbălat's Lemma (Lemma 4.2 on page 192 in [8]) or Barbashin-Krasovskii-LaSalle's theorem (Theorem 3.4 on page 115 in [8]). However, these tools cannot guarantee uniform attractivity and therefore these tools cannot be used for the verification of UAOS.

The counterexample in [14] showed that the existence of functions $V \in C^1(\Omega; \mathfrak{R}_+)$ with $V(0) = 0$, $a, b, c \in K_\infty$ for which

$$a(|h(x)|) \leq V(x) \leq b(|x|), \text{ for all } x \in \Omega \qquad (2.3)$$

$$\nabla V(x)f(x) \leq -c(|h(x)|), \text{ for all } x \in \Omega \tag{2.4}$$

does not guarantee UAOS on $\Omega$ even if in addition we assume that $h \in C^1(\Omega; \Re^k)$ and there exists a constant $M > 0$ with $|\nabla h(x)f(x)| \leq M$ for all $x \in \Omega$.

On the other hand, it is known that the existence of functions $V \in C^1(\Omega; \Re_+)$ with $V(0) = 0$, $a, b \in K_\infty$ and a continuous, positive definite function $\rho: \Re_+ \to \Re_+$ for which (2.3) holds and for which

$$\nabla V(x)f(x) \leq -\rho(V(x)), \text{ for all } x \in \Omega \tag{2.5}$$

does guarantee UAOS on $\Omega$ without any further assumption (see for instance [18,22] as well as Theorem 2.4 on page 90 in [6]). The key difference between (2.4) and (2.5) is that the dissipation rate in the former involves only the output norm, whereas the latter involves the Lyapunov function itself.

However, in practice it is difficult to obtain a dissipation rate, like in (2.5), that depends on the Lyapunov function itself. It is therefore reasonable to ask the following question:

"Assume that there exist functions $V \in C^1(\Omega; \Re_+)$, $W \in C^0(\Omega; \Re_+)$ and a continuous, positive definite function $\rho: \Re_+ \to \Re_+$ for which the following inequality holds:

$$\nabla V(x)f(x) \leq -\rho(W(x)), \text{ for all } x \in \Omega \tag{2.6}$$

What additional conditions on $V$ and $W$ can guarantee UAOS for (2.1), (2.2)?"

It is reasonable to expect that the dissipation rate will at least depend on the norm of the output, i.e., assume the existence of $a \in K_\infty$ such that

$$a(|h(x)|) \leq W(x), \text{ for all } x \in \Omega \tag{2.7}$$

Inequalities (2.6) and (2.7) guarantee that the dissipation rate is not zero when the output is non-zero (similarly to (2.4)). However, as recalled above, such an assumption is not sufficient for UOAS on $\Omega$.

2.3. Uniform Results

The following theorem provides sufficient Lyapunov-like conditions for UAOS of system (2.1), (2.2) without requiring a dissipation rate that involves the Lyapunov function and shows the missing link between AOS and UAOS in terms of the dissipation rate. Its proof is provided in Section 5.4.

**Theorem 1 (Lyapunov conditions for UAOS):** *Let $\Omega \subseteq S$ be a positively invariant set for (2.1) with $0 \in \Omega$. Suppose that there exist functions $a \in K_\infty$, $V, W \in C^1(\Omega; \Re_+)$ with $\sup\{V(x) + W(x) : x \in \Omega, |x| \leq s\} < +\infty$ for all $s \geq 0$, and a continuous, positive definite function $\rho: \Re_+ \to \Re_+$ such that inequalities (2.6), (2.7) hold as well as the following inequality:*

$$\nabla W(x)f(x) \leq 0, \text{ for all } x \in \Omega \tag{2.8}$$

*Then system (2.1), (2.2) is UAOS on $\Omega$.*

The above statement shows that a sufficient additional requirement in order for (2.6) and (2.7) to ensure UAOS is that the dissipation rate $W$ does not increase along the system's solutions.

**Remark 2.1:** Theorem 1 does not assume that (2.3) holds. Moreover, it is not assumed in Theorem 1 that $V(0) = 0$. The condition $\sup\{V(x)+W(x): x \in \Omega, |x| \leq s\} < +\infty$ for all $s \geq 0$ is automatically satisfied when $\Omega \subseteq \Re^n$ is a closed set. Notice that if (2.3) holds then the condition $\sup\{V(x): x \in \Omega, |x| \leq s\} < +\infty$ holds for all $s \geq 0$. Finally, notice that if (2.3) and (2.5) hold then all assumptions of Theorem 1 are valid with $W = V$, in which case we recover a classical sufficient condition for UAOS [6,18,22].

The use of Theorem 1 is illustrated in Section 3 by studying the stability properties of systems under adaptive control.
   The following theorem provides sufficient Lyapunov-like conditions for Lagrange output stability and Lyapunov output stability of system (2.1), (2.2).

**Proposition 1 (Lyapunov conditions for Lagrange and Lyapunov output stability):** *Let $\Omega \subseteq S$ be a positively invariant set for (2.1) with $0 \in \Omega$. Suppose that there exist functions $a \in K_\infty$, $W \in C^1(\Omega; \Re_+)$ with $W(0) = 0$, $\sup\{W(x): x \in \Omega, |x| \leq s\} < +\infty$ for all $s \geq 0$, such that inequalities (2.7), (2.8) hold. Then system (2.1), (2.2) is Lagrange output stable on $\Omega$ and Lyapunov output stable on $\Omega$.*

Proposition 1 is nothing new; it is a classical Lyapunov result. The only new thing in Proposition 1 is the use of the positively invariant set $\Omega \subseteq S$. However, we have included Proposition 1 in the paper for completeness and its proof is given in Section 5.3.

2.4. Non-Uniform Results

The main requirement in the previous section is (2.8), which imposes the requirement that $W$ does not increase along solutions. In case this condition is not satisfied, it is still possible to establish AOS using Barbălat's Lemma, which states that any integrable uniformly continuous function converges to zero. It has been observed in [20] that uniform continuity can be established under a uniform local integrability assumption. Here we show that uniform continuity is actually not required to conclude convergence to zero and that this requirement can be advantageously replaced by the following (less conservative) regularity property.

**Definition 2 (Quasi-uniform continuity):** *A function $f: \Re_+ \to \Re$ is said to be quasi-uniformly continuous if for every $\varepsilon > 0$ there exists $\delta(\varepsilon) > 0$ such that $f(t) - f(t_0) < \varepsilon$ for all $t \geq t_0 \geq 0$ with $t \leq t_0 + \delta(\varepsilon)$.*

Notice that any uniformly continuous function is quasi-uniformly continuous. Loosely speaking one could say that a function is quasi-uniformly continuous if it is "uniformly upper semi-continuous from the right" and "uniformly lower semi-continuous from the left". It is also worth stressing that any quasi-uniformly continuous function is Lebesgue measurable (see Fact 1 in [3]). Quasi-uniform continuity can be easily shown in practice using the following.

**Proposition 2 (Sufficient condition for quasi-uniform continuity):** *If for a function $f : \Re_+ \to \Re$ there exists a constant $M \geq 0$ such that the function $g : \Re_+ \to \Re$ defined by $g(t) = f(t) - M t$ for all $t \geq 0$ is non-increasing then $f : \Re_+ \to \Re$ is quasi-uniformly continuous. In particular, if $f \in C^1(\Re_+)$ and there exists $M > 0$ such that $\dot{f}(t) \leq M$ for all $t \geq 0$, then $f$ is quasi-uniformly continuous.*

Proposition 2 highlights the difference between uniform continuity and quasi-uniform continuity. Indeed, a sufficient condition for uniform continuity is that $|\dot{f}(t)| \leq M$ for all $t \geq 0$, which is commonly used in the application of Barbălat's lemma. As we will see below, requiring $\dot{f}(t)$ to be upper bounded (rather than both upper and lower bounded) sometimes proves useful. Proposition 2 is an immediate consequence of the definition of quasi-uniform continuity and its proof is therefore omitted.

We are now in position to present the following extension of Barbălat's Lemma, which does not require uniform continuity. Its proof is provided in Section 5.5.

**Lemma 3 (Extension of Barbălat's Lemma):** *Let $\rho : \Re_+ \to \Re_+$ be a continuous, positive definite function and let $f : \Re_+ \to \Re_+$ be a given function for which one of the following regularity requirements holds: either $f$ is quasi-uniformly continuous or $-f$ is quasi-uniformly continuous. Furthermore, suppose that $\int_0^{+\infty} \rho(f(t)) dt < +\infty$. Finally, suppose that either $\rho : \Re_+ \to \Re_+$ is non-decreasing or that $f : \Re_+ \to \Re_+$ is bounded. Then $\lim_{t \to +\infty} (f(t)) = 0$.*

Lemma 3 is less conservative than Barbălat's Lemma. Indeed, under the assumption that $\int_0^{+\infty} \rho(f(t)) dt < +\infty$, Barbălat's Lemma would allow us to conclude that $\lim_{t \to +\infty}(\rho(f(t))) = 0$ only if the mapping $t \to \rho(f(t))$ were uniformly continuous. On the other hand, Lemma 3 allows us to conclude that $\lim_{t \to +\infty}(f(t)) = 0$ (a stronger conclusion) if one of the functions $f$ or $-f$ is merely quasi-uniformly continuous.

Based on this extension, we have the following relaxed condition for AOS.

**Theorem 2 (Lyapunov conditions for AOS):** *Let $\Omega \subseteq S$ be a positively invariant set for (2.1) with $0 \in \Omega$. Suppose that there exist functions $a, b \in K_\infty$, $V, W \in C^1(\Omega; \Re_+)$ and a continuous, positive definite function $\rho : \Re_+ \to \Re_+$ for which inequalities (2.3), (2.6), (2.7) hold. Moreover, suppose that there exists a continuous function $\gamma : \Re_+ \to \Re_+$ such that one of the following inequalities hold:*

$$\nabla W(x) f(x) \leq \gamma(V(x)), \text{ for all } x \in \Omega \qquad (2.9)$$

*or*

$$\nabla W(x) f(x) \geq -\gamma(V(x)), \text{ for all } x \in \Omega \qquad (2.10)$$

*Finally, suppose that either $\rho : \Re_+ \to \Re_+$ is non-decreasing or that there exists a continuous function $\zeta : \Re_+ \to \Re_+$ such that the following inequality holds:*

$$W(x) \leq \zeta(V(x)), \text{ for all } x \in \Omega \qquad (2.11)$$

*Then system (2.1), (2.2) is AOS on $\Omega$.*

The proof of Theorem 2 is provided in Section 5.6.

**Remark 2.2: (i)** The proof of Theorem 2 uses Lemma 3 in an instrumental way. Notice that the assumptions of Theorem 2 guarantee that $t \to V(\phi(t, x_0))$ is bounded for every $x_0 \in \Omega$ and therefore that for every $x_0 \in \Omega$ there exists a constant $M(x_0) > 0$ such that either $\frac{d}{dt}W(\phi(t, x_0)) \leq M(x_0)$ holds for all $t \geq 0$ or $\frac{d}{dt}W(\phi(t, x_0)) \geq -M(x_0)$ holds for all $t \geq 0$. By virtue of Proposition 2 either the function $t \to W(\phi(t, x_0))$ is quasi-uniformly continuous or the function $t \to -W(\phi(t, x_0))$ is quasi-uniformly continuous. No bound on $\left|\frac{d}{dt}W(\phi(t, x_0))\right|$ is requested in the above statement, thus allowing to consider cases when $t \to W(\phi(t, x_0))$ is not uniformly continuous. This observation allows us to highlight a key difference between Theorems 1 and 2: the key assumption for UAOS is that $W$ does not increase along solutions, whereas AOS requires that it does not increase or decrease more than linearly in time.

**(ii)** If one of the following conditions holds: **(a)** $\Omega$ is a compact set, or **(b)** $V \in C^1(\Omega; \Re_+)$ is a radially unbounded function, then there exist continuous functions $\gamma, \zeta : \Re_+ \to \Re_+$ such that (2.9), (2.10) and (2.11) hold. This observation shows why Theorem 2 is less demanding than Barbashin-Krasovskii-LaSalle's theorem.

The following example illustrates the use of Theorem 2.

**Example 1 (Verification of AOS):** Consider the following 3-dimensional system

$$\begin{aligned}
\dot{y} &= -\left(1 + w^2\right) y + \frac{2zg(z, w)}{\left(1 + z^2\right)^2} \\
\dot{z} &= -g(z, w) y \\
\dot{w} &= w + |y| \\
x &= (y, z, w) \in \Re^3
\end{aligned} \qquad (2.12)$$

where $g : \Re^2 \to \Re$ is a bounded locally Lipschitz function. System (2.12) on $\Omega = \Re^3$ is forward complete. This can be shown by using the function $U(x) = \frac{1}{2}y^2 + \frac{1}{2}z^2 + \frac{1}{2}w^2$ which satisfies the following inequalities for all $x = (y, z, w) \in \Re^3$

$$\nabla U(x)\dot{x} = -\left(1+w^2\right)y^2 + \frac{2zg(z,w)}{\left(1+z^2\right)^2}y - zg(z,w)y + w^2 + w|y|$$

$$\leq -y^2 + 3|g(z,w)||z||y| + w^2 + |w||y| \tag{2.13}$$

$$\leq 9R^2z^2 + \frac{3}{2}w^2 \leq \left(18R^2+3\right)U(x)$$

where $R := \sup\left\{|g(z,w)|:(z,w) \in \Re^2\right\}$. It follows from (2.13) and Theorem 2 in [1] that system (2.12) on $\Omega = \Re^3$ is forward complete. Next we define the function $V(x) = \frac{1}{2}y^2 + \frac{z^2}{1+z^2}$, which satisfies (2.3) with $a(s) = \frac{1}{2}s^2$ and $b(s) = s^2$ for all $s \geq 0$. Moreover, $V$ satisfies the following inequality for all $x = (y,z,w) \in \Re^3$

$$\nabla V(x)\dot{x} = -\left(1+w^2\right)y^2 \leq -y^2 \tag{2.14}$$

Inequality (2.14) shows that (2.6) holds with $\rho(s) = 2s$ for all $s \geq 0$ (a non-decreasing function) and $W(x) = \frac{1}{2}y^2$. Furthermore, (2.7) holds with $a(s) = \frac{1}{2}s^2$ for all $s \geq 0$. Finally, we show that (2.9) holds. Indeed, by using the fact that $R := \sup\left\{|g(z,w)|:(z,w) \in \Re^2\right\}$, we obtain for all $x = (y,z,w) \in \Re^3$

$$\nabla W(x)\dot{x} = -\left(1+w^2\right)y^2 + \frac{2zg(z,w)}{\left(1+z^2\right)^2}y \leq -y^2 + R|y| \leq \frac{R^2}{4} \tag{2.15}$$

Inequality (2.15) shows that (2.9) holds with $\gamma(s) \equiv R^2/4$. Notice that since $\rho$ is non-decreasing, we do not have to show (2.11) (although it happens to hold with $\zeta(s) = s$ for all $s \geq 0$). Therefore Theorem 2 guarantees that system (2.12) is AOS on $\Omega = \Re^3$.

Barbashin-Krasovskii-LaSalle's theorem could not have been used in order to prove AOS on $\Omega = \Re^3$ for (2.11), because it admits unbounded solutions. To see this notice that the differential inequality $\dot{w} \geq w$ holds for all $x = (y,z,w) \in \Re^3$. Thus we get $w(t) \geq \exp(t)w(0)$ for all $t \geq 0$, and consequently, when $w(0) > 0$ the solutions of (2.11) are unbounded.

In the same way, Barbălat's Lemma could not have been easily used in order to prove AOS on $\Omega = \Re^3$ for (2.11), because (2.14) shows that $\dot{W}(x)$ is bounded from above but not necessarily bounded from below. In other words, we cannot (easily) prove that the mapping $t \to W(\phi(t,x))$ is uniformly continuous.  ◁

# 3. Application to Adaptive Control

Adaptive control is a particularly relevant situation in which AOS is useful. Such control strategies rely on a state extension to provide a dynamical estimate of unknown parameters. In this context, it is not always necessary to precisely estimate these unknown parameters to get a satisfactory behavior of the system's solutions. In other words, only part of the state variables of the closed-loop system are requested to converge to zero, which can be rephrased as an AOS objective. In most adaptive control literature (e.g. [9]), convergence of this state variables is concluded through Barbălat's lemma or Barbashin-Krasovskii-LaSalle arguments, thus impeding to show UAOS and generating potential robustness issues. We show below how Theorem 1 can be employed to guarantee UAOS under adaptive control. To that aim, consider the following finite-dimensional control system

$$\dot{y} = f(y) + g(y)u + g(y)(\varphi(y))'\theta \tag{3.1}$$
$$y \in \Re^n, \theta \in \Re^p, u \in \Re$$

where $f, g : \Re^n \to \Re^n$, $\varphi : \Re^n \to \Re^p$ are locally Lipschitz vector fields with $f(0) = 0$ and $\varphi(0) = 0$, $y \in \Re^n$ denotes the state, $u \in \Re$ is the control input and $\theta \in \Re^p$ is the constant vector of unknown parameters. The control system (3.1) satisfies the "matching condition" (following the terminology in [9]) since the uncertain parameters are in the span of the control.

The following assumption is used for the control system (3.1).

**(H)** *There exist functions $P \in C^2(\Re^n; \Re_+)$ being positive definite and radially unbounded, $Q \in C^1(\Re^n; \Re_+)$ being positive definite and a locally Lipschitz function $k \in C^0(\Re^n; \Re)$ with $k(0) = 0$ such that*

$$\nabla P(y) f(y) + \nabla P(y) g(y) k(y) \leq -Q(y), \text{ for all } y \in \Re^n \tag{3.2}$$

When Assumption (H) holds then it is possible to regulate the state to zero by applying a standard adaptive control scheme. To see this, consider the adaptive controller

$$u = k(y) - (\varphi(y))'\hat{\theta}$$
$$\frac{d\hat{\theta}}{dt} = \gamma^{-1}\nabla P(y) g(y) \varphi(y) \tag{3.3}$$

where $\gamma > 0$ is a constant parameter of the controller. System (3.1) in closed loop with (3.3) then reads:

$$\dot{y} = f(y) + g(y)\left(k(y) - (\varphi(y))' z\right)$$
$$\dot{z} = \gamma^{-1}\nabla P(y) g(y) \varphi(y) \tag{3.4}$$
$$x = (y, z) \in \Re^n \times \Re^p$$

with $z = \hat{\theta} - \theta$ and output map defined by

$$h(x) = y \tag{3.5}$$

The closed-loop system (3.4)-(3.5) is AOS on $\Re^n \times \Re^p$. The proof of this fact follows the analysis given in [9] using the Lyapunov function $V(x) := P(y) + \frac{\gamma}{2}|z|^2$, which (by virtue of (3.2), (3.4)) satisfies the following differential inequality

$$\frac{d}{dt}V(x) \leq -Q(y), \text{ for all } x = (y,z) \in \Re^n \times \Re^p \tag{3.6}$$

However, the convergence rate of $y \in \Re^n$ to zero is not necessarily uniform with respect to the initial condition, even considering initial states in a bounded neighborhood of the origin.

In order to be able to guarantee UAOS on a prescribed set $\Omega \subseteq \Re^n \times \Re^p$ we need an additional assumption.

**(A)** *There exists a locally Lipschitz, positive function $\mu \in C^0\left(\Re^n;(0,+\infty)\right)$ such that the following inequality holds:*

$$|\varphi(y)|^2 \leq \mu(y)Q(y), \text{ for all } y \in \Re^n \tag{3.7}$$

Assumption (A) is automatically guaranteed if there exists a constant $q > 0$ such that the inequality $Q(y) \geq q|y|^2$ holds for all $y \in \Re^n$ in a neighborhood of $0 \in \Re^n$. Both Assumptions (H) and (A) are automatically guaranteed if the feedback law $u = k(y)$ achieves global asymptotic stabilization and local exponential stabilization of $0 \in \Re^n$ for the control system $\dot{y} = f(y) + g(y)u$ (recall Proposition 4.4 in [7]).

When Assumptions (H) and (A) hold then it is possible to slightly modify the adaptive control scheme (3.3) in the following way:

$$u = k(y) - (\varphi(y))' \hat{\theta} - L\mu(y)\nabla P(y)g(y)$$
$$\frac{d\hat{\theta}}{dt} = \gamma^{-1}\nabla P(y)g(y)\varphi(y) \tag{3.8}$$

where $\gamma, L > 0$ are constant parameters of the controller. With this new adaptive strategy, the system (3.1) in closed loop with (3.8) reads:

$$\dot{y} = f(y) + g(y)\left(k(y) - L\mu(y)\nabla P(y)g(y) - (\varphi(y))' z\right)$$
$$\dot{z} = \gamma^{-1}\nabla P(y)g(y)\varphi(y) \tag{3.9}$$
$$x = (y,z) \in \Re^n \times \Re^p$$

with $z = \hat{\theta} - \theta$ and output map defined by (3.5). The result below states that the closed-loop system (3.9) is UAOS on the set

$$\Omega := \left\{ x = (y,z) \in \Re^n \times \Re^p : V(x) = P(y) + \frac{\gamma}{2}|z|^2 \leq \gamma L \right\} \tag{3.10}$$

Notice that the size of the set $\Omega \subseteq \Re^n \times \Re^p$ is determined by the controller parameters and can be made as large as desired (by picking sufficiently large constants $\gamma, L > 0$).

**Theorem 3 (Semiglobal UAOS through adaptive control):** *Suppose that Assumptions (H) and (A) hold. Then for every $\gamma, L > 0$ system (3.9) with output (3.5) is AOS on $\Re^n \times \Re^p$ and UAOS on $\Omega \subseteq \Re^n \times \Re^p$, where $\Omega$ is defined by (3.10).*

The proof of Theorem 3 is provided in Section 5.7. When Assumptions (H) and (A) hold and the vector of unknown parameters $\theta \in \mathfrak{R}^p$ belongs to a *known bounded set* $\Theta \subseteq \mathfrak{R}^p$ then it is possible to achieve robust global asymptotic stabilization (even when $\theta \in \mathfrak{R}^p$ is time-varying; see [6]) of the equilibrium point $y = 0$ of system (3.1) by using the feedback law

$$u = k(y) - L\mu(y)\nabla P(y)g(y)$$

where $L > 0$ is a sufficiently large constant (that depends on the "size" of the set $\Theta \subseteq \mathfrak{R}^p$). It is clear that the adaptive controller (3.8) is a combination of the above feedback law and the adaptive controller (3.3). When $L \to 0$ then the controller (3.8) "tends" to the adaptive controller (3.3).

Theorem 3 is important because Theorem 1 in conjunction with Lemma 1 and Lemma 2 guarantees that there exist functions $\sigma \in KL$, $\zeta \in K_\infty$ such that the following estimates hold for the solution closed-loop system (3.1) with (3.8):

$$|y(t)| \leq \zeta\left(|y(0)| + |\hat{\theta}(0) - \theta|\right), \text{ for all } t \geq 0, \ (y(0), \hat{\theta}(0)) \in \mathfrak{R}^n \times \mathfrak{R}^p \text{ and } \theta \in \mathfrak{R}^p \quad (3.11)$$

$$\lim_{t \to +\infty}\left(|y(t)|\right) = 0, \text{ for all } (y(0), \hat{\theta}(0)) \in \mathfrak{R}^n \times \mathfrak{R}^p \text{ and } \theta \in \mathfrak{R}^p \quad (3.12)$$

$$|y(t)| \leq \sigma\left(|y(0)| + |\hat{\theta}(0) - \theta|, t\right),$$

for all $t \geq 0$, $\theta \in \mathfrak{R}^p$ and $(y(0), \hat{\theta}(0)) \in \mathfrak{R}^n \times \mathfrak{R}^p$ with $(y(0), \hat{\theta}(0) - \theta) \in \Omega$ \quad (3.13)

It should be noticed that the functions $\sigma \in KL$, $\zeta \in K_\infty$ are independent of $\theta \in \mathfrak{R}^p$ (since system (3.9) does not depend on $\theta \in \mathfrak{R}^p$), thus guaranteeing the same qualitative behavior of solutions regardless of the actual value of the parameters. Moreover, it should be noticed that estimate (3.13) can be used in conjunction with Lemma 9 in [14] in order to establish an Input-to-Output Stability (IOS) estimate for the closed-loop system (3.1) with (3.8) under the effect of possible modeling errors. Therefore, the "Lyapunov redesign" procedure described above enhances the robustness properties of the closed-loop system. Finally, it should be noticed that a *KL* estimate like (3.13) is rarely available in the literature for systems under adaptive control, in which no uniform convergence to the origin is usually ensured.

## 4. Output Stability for Abstract Systems

We now investigate on how the results of Section 2, developed in a finite-dimensional context, can be extended to more general classes of dynamical systems.

### 4.1. Definitions

We start by providing the definition of an abstract, continuous-time, deterministic, autonomous, dynamical system with output.

**Definition 3 (Abstract dynamical system):** *Let $X, Y$ be normed linear spaces with norms $\|\cdot\|_X$, $\|\cdot\|_Y$, respectively. A (continuous-time, deterministic, autonomous) dynamical system with output is a triplet $\Sigma = (S, h, \phi)$ that consists of a non-empty set $S \subseteq X$ (the state space), a mapping $h: S \to Y$ (the output map) and a mapping $\phi: \mathfrak{R}_+ \times S \to S$ that satisfies: (i) the identity property, i.e., $\phi(0, x) = x$ for all $x \in S$, and (ii) the semigroup property, i.e., $\phi(t, \phi(s, x)) = \phi(t + s, x)$ for all $x \in S$ and $t, s \geq 0$.*

Definition 3 is a specialization of the more general definitions of abstract control systems given in [6,12] and [16]. It can be used for the study of finite-dimensional systems (when $X = \Re^n$), delay systems (where $X = C^0\left([-r,0];\Re^n\right)$) or systems described by Partial Differential Equations (where $X$ is an appropriate functional space; see [12]).

Throughout this section, we assume the following.

**Standing Assumption 1:** $\Omega \subseteq S$ *is a non-empty set for which* $\phi(t,\Omega) \subseteq \Omega$ *for all* $t \geq 0$. *Moreover,* $0 \in \Omega$, $h(0) = 0$ *and* $\phi(t,0) = 0$ *for all* $t \geq 0$.

The notions of output stability can be extended to such general dynamical systems with outputs. In what follows, we use the notation $y(t,x_0) = h(\phi(t,x_0))$ for all $t \geq 0$, $x_0 \in S$ and $B_R = \left\{ x \in X : \|x\|_X < R \right\}$ for all $R > 0$.

**Definition 4 (Output stability properties):** *Let* $\Sigma = (S,h,\phi)$ *be a dynamical system with output and let* $\Omega \subseteq S$ *be a set satisfying Standing Assumption 1. We say that system* $\Sigma$ *is*
  i) *Lagrange output stable on* $\Omega$ *if for every* $R > 0$ *the set* $\left\{ \|y(t,x_0)\|_Y : x_0 \in B_R \cap \Omega, t \geq 0 \right\}$ *is bounded.*
  ii) *Lyapunov output stable on* $\Omega$ *if for every* $\varepsilon > 0$ *there exists* $\delta(\varepsilon) > 0$ *such that for all* $x_0 \in B_{\delta(\varepsilon)} \cap \Omega$ *it holds that* $\|y(t,x_0)\|_Y \leq \varepsilon$ *for all* $t \geq 0$.
  iii) *Asymptotically Output Stable (AOS) on* $\Omega$ *if system (2.1), (2.2) is Lagrange and Lyapunov output stable on* $\Omega$ *and* $\lim_{t \to +\infty}\left(y(t,x_0)\right) = 0$ *for all* $x_0 \in \Omega$.
  iv) *Uniformly Asymptotically Output Stable (UAOS) on* $\Omega$ *if system (2.1), (2.2) is Lagrange and Lyapunov output stable on* $\Omega$ *and for every* $\varepsilon, R > 0$ *there exists* $T(\varepsilon,R) > 0$ *such that for all* $x_0 \in B_R \cap \Omega$ *it holds that* $\|y(t,x_0)\|_Y \leq \varepsilon$ *for all* $t \geq T(\varepsilon,R)$.

The following two lemmas are similar to Lemmas 1 and 2 and provide characterizations of the stability notions of Definition 3 in terms of comparison functions. These two results are respectively established in Sections 5.1 and 5.2.

**Lemma 4 (Output stability through $K_\infty$ estimate):** *Let* $\Sigma = (S,h,\phi)$ *be a dynamical system with output and let* $\Omega \subseteq S$ *be a set satisfying Standing Assumption 1. System* $\Sigma$ *is Lagrange output stable on* $\Omega$ *and Lyapunov output stable on* $\Omega$ *if and only if there exists a function* $\zeta \in K_\infty$ *such that the following estimate holds for all* $x_0 \in \Omega$ *and* $t \geq 0$:

$$\|y(t,x_0)\|_Y \leq \zeta\left(\|x_0\|_X\right)$$

**Lemma 5 (UAOS through $KL$ estimate):** *Let* $\Sigma = (S,h,\phi)$ *be a dynamical system with output and let* $\Omega \subseteq S$ *be a set satisfying Standing Assumption 1. System* $\Sigma$ *is UAOS on* $\Omega$ *if and only if there exists a function* $\sigma \in KL$ *such that the following estimate holds for all* $x_0 \in \Omega$ *and* $t \geq 0$:

$$\|y(t,x_0)\|_Y \leq \sigma\left(\|x_0\|_X, t\right)$$

### 4.2. Uniform Results

The following theorem is similar to Theorem 1 and provides sufficient conditions for UAOS of an abstract dynamical system $\Sigma$ with output. The proof of Theorem 4 is exactly the same as the proof of Theorem 1 and is omitted.

**Theorem 4 (Lyapunov conditions for UAOS):** *Let $\Sigma = (S, h, \phi)$ be a dynamical system with output and let $\Omega \subseteq S$ be a set satisfying Standing Assumption 1. Suppose that there exist a function $a \in K_\infty$, continuous functionals $V, W : \Omega \to \Re_+$ with $\sup\{V(x) + W(x) : x \in \Omega, \|x\|_X \leq s\} < +\infty$ for all $s \geq 0$, and a continuous, positive definite function $\rho : \Re_+ \to \Re_+$ such that*

$$a(\|h(x)\|_Y) \leq W(x), \text{ for all } x \in \Omega \tag{4.1}$$

$$V(\phi(t,x)) \leq V(x) - \int_0^t \rho(W(\phi(s,x)))\,ds, \text{ for all } x \in \Omega, t \geq 0 \tag{4.2}$$

*Moreover, suppose that for every $x \in \Omega$ the mapping $t \to W(\phi(t,x))$ is non-increasing. Then system $\Sigma$ is UAOS on $\Omega$.*

The following statement is the analogue of Proposition 1. Its proof is similar to the proof of Proposition 1 and is omitted.

**Proposition 3 (Lyapunov conditions for Lagrange and Lyapunov output stability):** *Let $\Sigma = (S, h, \phi)$ be a dynamical system with output and let $\Omega \subseteq S$ be a set satisfying Standing Assumption 1. Suppose that there exist functions $a, b \in K_\infty$ and a functional $V : \Omega \to \Re_+$ for which the following property holds:*

$$a(\|h(x)\|_Y) \leq V(x) \leq b(\|x\|_X), \text{ for all } x \in \Omega \tag{4.3}$$

*Moreover, suppose that for every $x \in \Omega$ the mapping $t \to V(\phi(t,x))$ is non-increasing. Then system $\Sigma$ is Lagrange output stable on $\Omega$ and Lyapunov output stable on $\Omega$.*

### 4.3. Non-Uniform Results

The following theorem is similar to Theorem 2 and provides sufficient conditions for AOS of an abstract dynamical system $\Sigma$ with output. The proof of Theorem 5 is exactly the same as the proof of Theorem 2 and is omitted.

**Theorem 5 (Lyapunov conditions for AOS):** *Let $\Sigma = (S, h, \phi)$ be a dynamical system with output and let $\Omega \subseteq S$ be a set satisfying Standing Assumption 1. Suppose that there exist a continuous, positive definite function $\rho : \Re_+ \to \Re_+$, functions $a, b \in K_\infty$ and functionals $V, W : \Omega \to \Re_+$ for which inequalities (4.1), (4.2), (4.3) hold. Moreover, suppose that either for every $x \in \Omega$ the mapping $t \to W(\phi(t,x))$ is quasi-uniformly continuous or that for every $x \in \Omega$ the mapping $t \to -W(\phi(t,x))$ is quasi-uniformly continuous. Finally, suppose that either $\rho : \Re_+ \to \Re_+$ is non-decreasing or that for every $x \in \Omega$ the mapping $t \to W(\phi(t,x))$ is bounded. Then system $\Sigma$ is AOS on $\Omega$.*

**Remark 3.1:** Lemma 6.3 on page 186 in [2] in conjunction with the semigroup property implies that a necessary and sufficient condition for the integral inequality (4.2) under the assumption that the mapping $t \to V(\phi(t,x))$ is lower semi-continuous and the mapping $t \to W(\phi(t,x))$ is continuous, is the following differential inequality:

$$\limsup_{t \to 0^+} \left( \frac{V(\phi(t,x)) - V(x)}{t} \right) \leq -\rho(W(x)), \text{ for all } x \in \Omega$$

which allows to test assumption (4.2) using the upper-right Dini derivative of $V$ along the system's solutions.

4.4. Output Stability for Delay Systems

We next show how to use the proposed result in the context of time-delay systems

$$\dot{x}(t) = f(x_t)$$
$$x(t) \in \mathfrak{R}^n \tag{4.4}$$

where $f: X = C^0\left([-r,0];\mathfrak{R}^n\right) \to \mathfrak{R}^n$, with $r > 0$ being a constant, is a Lipschitz mapping on bounded sets of $X$ with $f(0) = 0$ and $x_t \in X$ is the state defined by $(x_t)(s) = x(t+s)$ for all $s \in [-r,0]$. The constant $r$ can be picked as any constant greater than or equal to the maximal delay involved in the system dynamics. Let $h: X \to \mathfrak{R}^k$ be a continuous mapping with $h(0) = 0$. We consider system (4.4) with output given by the following equation:

$$y(t) = h(x_t) \tag{4.5}$$

The following corollaries are direct consequences of Theorem 3 and Theorem 4. In what follows for every $x \in X$, the norm of $\|x\|$ is given by $\|x\| = \max_{s \in [-r,0]} (|x(s)|)$.

**Corollary 1 (UAOS for time-delay systems):** *Let $\Omega \subseteq X$ with $0 \in \Omega$ be a positively invariant set for system (4.4), (4.5) and suppose that there exist a function $a \in K_\infty$, continuous functionals $V, W : \Omega \to \mathfrak{R}_+$ with $\sup\{V(x) + W(x) : x \in \Omega, \|x\| \leq s\} < +\infty$ for all $s \geq 0$, and a continuous, positive definite function $\rho: \mathfrak{R}_+ \to \mathfrak{R}_+$ such that the following inequalities hold for all $x_0 \in \Omega$:*

$$a(|h(x_0)|) \leq W(x_0) \tag{4.6}$$

$$\limsup_{t \to 0^+} \left( \frac{V(\phi(t,x_0)) - V(x_0)}{t} \right) \leq -\rho(W(x_0)) \tag{4.7}$$

$$\limsup_{t \to 0^+} \left( \frac{W(\phi(t,x_0)) - W(x_0)}{t} \right) \leq 0 \tag{4.8}$$

*where $\phi(t, x_0) \in \Omega$ denotes the solution $x_t \in \Omega$ of (4.5) with initial condition $x_0 \in \Omega$. Then system (4.4), (4.5) is UAOS on $\Omega$.*

It is worth stressing that the upper-right Dini derivatives appearing in (4.7) and (4.8) can be computed with no knowledge of the system's solutions, for instance using Driver's derivative: see for instance [15].

**Corollary 2 (AOS for time-delay systems):** *Let $\Omega \subseteq X$ with $0 \in \Omega$ be a positively invariant set for system (4.4), (4.5) and suppose that there exist a continuous, positive definite function $\rho : \Re_+ \to \Re_+$, a continuous function $\gamma : \Re_+ \to \Re_+$, functions $a,b \in K_\infty$ and functionals $V,W : X \to \Re_+$ for which inequalities (4.6), (4.7) hold as well as the following inequality:*

$$a(|h(x)|) \leq V(x) \leq b(\|x\|), \text{ for all } x \in \Omega \tag{4.9}$$

*Furthermore, assume that one of the following inequalities holds:*

$$\limsup_{t \to 0^+} \left( \frac{W(\phi(t,x_0)) - W(x_0)}{t} \right) \leq \gamma(V(x_0)), \text{ for all } x_0 \in \Omega \tag{4.10}$$

or

$$\limsup_{t \to 0^+} \left( \frac{W(x_0) - W(\phi(t,x_0))}{t} \right) \leq \gamma(V(x_0)), \text{ for all } x_0 \in \Omega \tag{4.11}$$

*where $\phi(t,x_0) \in \Omega$ denotes the solution $x_t \in \Omega$ of (4.5) with initial condition $x_0 \in \Omega$. Finally, suppose that either $\rho : \Re_+ \to \Re_+$ is non-decreasing or that there exists a continuous function $\zeta : \Re_+ \to \Re_+$ such that the following inequality holds:*

$$W(x) \leq \zeta(V(x)), \text{ for all } x \in \Omega \tag{4.12}$$

*Then system (4.4), (4.5) is AOS on $\Omega$.*

The following example illustrates the use of Corollary 1.

**Example 2 (Verification of UAOS):** Consider the time-delay system

$$\begin{aligned}
\dot{x}_1(t) &= -px_1(t) + qx_1(t-r) - g(x(t))x_1(t)x_2(t) \\
\dot{z}(t) &= g(x(t))x_1^2(t) \\
x(t) &= (x_1(t), x_2(t)) \in \Re^2
\end{aligned} \tag{4.13}$$

where $r > 0$, $p > 0$, $q \in \Re$ are constants and $g : \Re^2 \to \Re$ is a locally Lipschitz function. We assume that there exist constants $R, \sigma, Q > 0$ such that the following inequalities hold:

$$\frac{4\lambda(p-Q-\lambda)^2 + \sigma^2 Q}{2\sigma(p-Q-\lambda)} \leq p + g(x)z, \text{ for all } x \in \Re^2 \text{ with } |x| \leq R \tag{4.14}$$

where

$$\lambda := \frac{q^2 \exp(\sigma r)}{4Q} < p - Q \tag{4.15}$$

Consider the continuous functionals $V : C^0([-r,0]; \Re^2) = X \to \Re_+$, $W : X \to \Re_+$ defined by the following equations for all $x \in X$:

$$V(x) := \frac{1}{2}x_1^2(0) + Q\int_{-r}^{0} \exp(\sigma\theta)x_1^2(\theta)d\theta + \frac{1}{2}x_2^2(0) \tag{4.16}$$

$$W(x) := \frac{1}{2}x_1^2(0) + K\int_{-r}^{0} \exp(\sigma\theta)x_1^2(\theta)d\theta \tag{4.17}$$

where

$$K := \frac{\sigma Q}{2(p-Q-\lambda)} \tag{4.18}$$

Moreover, define the set

$$\Omega := \left\{ x \in X : V(x) \le \frac{R^2}{2} \right\} \tag{4.19}$$

and the output map for all $x \in X$:

$$h(x) = x_1(0) \tag{4.20}$$

Definitions (4.16), (4.17) and assumption (4.15) imply that inequalities (4.6), (4.9) hold with $a(s) = s^2/2$ and $b(s) = \left(Q + \frac{1}{2}\right)s^2$ for $s \ge 0$. Notice that the derivative of $V$ along the solutions of (4.13) satisfies the following inequalities for all $x \in X$:

$$\begin{aligned}\limsup_{t\to 0^+}\left(\frac{V(\phi(t,x))-V(x)}{t}\right) &= -(p-Q)x_1^2(0) + qx_1(-r)x_1(0) \\ &\quad -\sigma Q \int_{-r}^{0} \exp(\sigma s)x_1^2(s)ds - Q\exp(-\sigma r)x_1^2(-r) \\ &\le -2(p-Q-\lambda)W(x)\end{aligned} \tag{4.21}$$

For the derivation of (4.21) we have used definitions (4.15), (4.16), (4.17), (4.18) and the fact that $qx_1(-r)x_1(0) \le \frac{q^2}{4Q\exp(-\sigma r)}x_1^2(0) + Q\exp(-\sigma r)x_1^2(-r)$ for all $x_1(0), x_1(-r) \in \Re$. Inequality (4.21) and standard arguments allow us to conclude that $\Omega \subseteq X$ is a positively invariant set for system (4.13). Moreover, inequality (4.21) shows that (4.7) holds with $\rho(s) = 2\left(p - Q - \frac{q^2}{4Q\exp(-\sigma r)}\right)s$ for $s \ge 0$.

Finally, we notice that the derivative of $W$ along the solutions of (4.13) satisfies the following inequalities for all $x \in \Omega$:

$$\limsup_{t \to 0^+} \left( \frac{W(\phi(t,x)) - W(x)}{t} \right) = -\left(p - K + g(x(0))x_2(0)\right)x_1^2(0)$$
$$+ qx_1(-r)x_1(0) - \sigma K \int_{-r}^{0} \exp(\sigma s) x_1^2(s) ds - K \exp(-\sigma r) x_1^2(-r) \le 0 \tag{4.22}$$

For the derivation of (4.22) we have used definitions (4.17), (4.18), the fact that $qx_1(-r)x_1(0) \le K \exp(-\sigma r) x_1^2(-r) + \frac{q^2}{4K \exp(-\sigma r)} x_1^2(0)$ for all $x_1(0), x_1(-r) \in \Re$, inequality (4.14) and the fact that for all $x \in \Omega$ it holds that $|x(0)| \le R$ (a direct consequence of definitions (4.16), (4.19)).

It follows from Corollary 1 that system (4.13) is UAOS on $\Omega$. ◁

## 5. Proofs

### 5.1. Proof of Lemma 4

If $\Sigma$ is Lagrange output stable on $\Omega$ and Lyapunov output stable on $\Omega$ then we define $a(s) := \sup\{\|y(t,x)\|_Y : t \ge 0, x \in \Omega, \|x\|_X \le s\}$ for all $s \ge 0$. The function $a : \Re_+ \to \Re_+$ is well-defined by virtue of Lagrange output stability on $\Omega$ with $a(0) = 0$. Moreover, $a : \Re_+ \to \Re_+$ is continuous at $s = 0$ by virtue of Lyapunov output stability on $\Omega$. Finally, we notice that $a : \Re_+ \to \Re_+$ is non-decreasing. The existence of a function $\zeta \in K_\infty$ such that the estimate $\|y(t,x)\|_Y \le \zeta(\|x\|_X)$ holds for all $x \in \Omega$ and $t \ge 0$ is a direct consequence of Lemma 2.4 on page 65 in [6].

If there exists a function $\zeta \in K_\infty$ such that the estimate $\|y(t,x)\|_Y \le \zeta(\|x\|_X)$ holds for all $x \in \Omega$ and $t \ge 0$ then we can verify in a straightforward way Lagrange output stability on $\Omega$ and Lyapunov output stability on $\Omega$. The proof is complete. ◁

### 5.2. Proof of Lemma 5

If system $\Sigma$ is UAOS on $\Omega$ then we can define the functions $a : \Re_+ \to \Re_+$ and $M : \Re_+ \times \Re_+ \to \Re_+$ by means of the formulae

$$a(s) := \sup\{\|y(t,x)\|_Y : t \ge 0, x \in \Omega, \|x\|_X \le s\} \text{ for } s \ge 0$$

$$M(t,s) := \sup\{\|y(t,x)\|_Y : x \in \Omega, \|x\|_X \le s\} \text{ for } s \ge 0.$$

Working with $a : \Re_+ \to \Re_+$ and $M : \Re_+ \times \Re_+ \to \Re_+$, we can show the existence of a function $\sigma \in KL$ such that the estimate $\|y(t,x)\|_Y \le \sigma(\|x\|_X, t)$ holds for all $x \in \Omega$ and $t \ge 0$ by following exactly the same steps of the proof of Lemma 2.6 in [6].

If there exists a function $\sigma \in KL$ such that the estimate $\|y(t,x)\|_Y \le \sigma(\|x\|_X, t)$ holds for all $x \in \Omega$ and $t \ge 0$ then we can verify in a straightforward way UAOS on $\Omega$. The proof is complete. ◁

## 5.3. Proof of Proposition 1

There exists a function $b \in K_\infty$ such that

$$W(x) \leq b(|x|), \text{ for all } x \in \Omega \tag{5.1}$$

To see this, we may define the non-decreasing, non-negative function $\tilde{b}(s) := \sup\{W(x) : x \in \Omega, |x| \leq s\}$ for all $s \geq 0$. By continuity of $W$ and the fact that $W(0) = 0$, it follows that $\lim_{s \to 0^+}(\tilde{b}(s)) = \tilde{b}(0) = 0$. Moreover, the definition of $\tilde{b}(s)$ guarantees that $W(x) \leq \tilde{b}(|x|)$ for all $x \in \Omega$. The existence of $b \in K_\infty$ for which inequality (5.1) holds is a direct consequence of Lemma 2.4 on page 65 in [6].

Let (arbitrary) $x_0 \in \Omega$ be given. It follows from (2.8) that

$$\frac{d}{dt}W(\phi(t, x_0)) \leq 0, \text{ for all } t \geq 0 \tag{5.2}$$

Therefore, by virtue of (5.2) and since $\phi(0, x_0) = x_0$, the following estimate holds for all $t \geq 0$:

$$W(\phi(t, x_0)) \leq W(x_0) \tag{5.3}$$

Using (5.3) in conjunction with (2.7), (5.1), we obtain the estimate:

$$a(|y(t, x_0)|) \leq b(|x_0|) \tag{5.4}$$

By virtue of Lemma 1, estimate (5.4) shows Lagrange output stability on $\Omega$ and Lyapunov output stability on $\Omega$.  ◁

## 5.4. Proof of Theorem 1

Notice that inequality (2.6) for $x = 0$ implies that $W(0) = 0$. Therefore, Proposition 1 implies Lagrange output stability on $\Omega$ and Lyapunov output stability on $\Omega$.

Notice that (2.6) implies the following differential inequality for all $x_0 \in \Omega$:

$$\frac{d}{dt}V(\phi(t, x_0)) \leq -\rho(W(\phi(t, x_0))), \text{ for all } t \geq 0 \tag{5.5}$$

Moreover, (2.8) implies that the following estimate holds for all $x_0 \in \Omega$, $t \geq 0$ and $t_0 \in [0, t]$:

$$W(\phi(t, x_0)) \leq W(\phi(t_0, x_0)) \tag{5.6}$$

Let (arbitrary) $\varepsilon, R > 0$ be given. Define:

$$\tilde{\rho}(\varepsilon, R) := \min\{\rho(s) : a(\varepsilon) \leq s \leq a(\varepsilon) + \sup\{W(x) : x \in \Omega, |x| \leq R\}\} \tag{5.7}$$

$$T(\varepsilon, R) := \frac{1 + \sup\{V(x) : x \in \Omega, |x| \leq R\}}{\tilde{\rho}(\varepsilon, R)} \tag{5.8}$$

We next prove by contradiction that $\sup\{|y(t,x_0)|: x_0 \in B_R \cap \Omega, t \geq T(\varepsilon, R)\} \leq \varepsilon$, thus establishing UAOS. To that end, suppose that $\sup\{|y(t,x_0)|: x_0 \in B_R \cap \Omega, t \geq T(\varepsilon, R)\} > \varepsilon$. Then there exists $x_0 \in B_R \cap \Omega$ and $t \geq T(\varepsilon, R)$ such that $|y(t,x_0)| > \varepsilon$. By virtue of (2.7) we conclude that $W(\phi(t,x_0)) > a(\varepsilon)$. It follows from (5.6) that

$$a(\varepsilon) \leq W(\phi(s,x_0)) \leq W(x_0), \text{ for all } s \in [0,t] \tag{5.9}$$

Since $x_0 \in B_R \cap \Omega$, we get from (5.9) that

$$a(\varepsilon) \leq W(\phi(s,x_0)) \leq \sup\{W(x): x \in \Omega, |x| \leq R\}, \text{ for all } s \in [0,t] \tag{5.10}$$

and consequently by using definition (5.7)

$$\tilde{\rho}(\varepsilon, R) \leq \rho(W(\phi(s,x_0))), \text{ for all } s \in [0,t] \tag{5.11}$$

Using (5.5) we obtain:

$$V(\phi(t,x_0)) \leq V(x_0) - \int_0^t \rho(W(\phi(s,x_0))) ds \tag{5.12}$$

Combining (5.11) and (5.12), we get the inequality

$$V(\phi(t,x_0)) \leq V(x_0) - t\,\tilde{\rho}(\varepsilon, R) \tag{5.13}$$

Since $x_0 \in B_R \cap \Omega$, we get from (5.13)

$$V(\phi(t,x_0)) \leq \sup\{V(x): x \in \Omega, |x| \leq R\} - t\,\tilde{\rho}(\varepsilon, R) \tag{5.14}$$

Inequality (5.14) in conjunction with the fact that $t \geq T(\varepsilon, R)$ and definition (5.8) gives

$$V(\phi(t,x_0)) < 0$$

which contradicts the fact that $V(\phi(t,x_0)) \geq 0$. ◁

5.5. Proof of Lemma 3

The proof is made by contradiction. Suppose that there exists $\varepsilon > 0$ and an increasing sequence $\{t_i \geq 0 : i = 1, 2, ...\}$ with $\lim_{i \to +\infty}(t_i) = +\infty$ such that $f(t_i) \geq \varepsilon$.

We next distinguish the following cases.

Case 1: $f$ is quasi-uniformly continuous

Since $f$ is a quasi-uniformly continuous function there exists $\delta > 0$ such that $f(t) - f(s) < \varepsilon/2$ for all $s \in [t-\delta, t]$.

Without loss of generality (by taking a subsequence if necessary) we may assume that $t_1 \geq 1+\delta$ and $t_{i+1} \geq t_i + \delta + 1$ for all $i = 1, 2, \ldots$. Notice that

$$\frac{\varepsilon}{2} \leq f(t_i) - \frac{\varepsilon}{2} < f(s) \text{ for all } s \in [t_i - \delta, t_i], \ i = 1, 2, \ldots. \tag{5.15}$$

If $\rho : \Re_+ \to \Re_+$ is non-decreasing then (5.15) implies $\rho\left(\frac{\varepsilon}{2}\right) \leq \rho(f(s))$ for all $s \in [t_i - \delta, t_i]$, $i = 1, 2, \ldots$. Since the intervals $[t_i - \delta, t_i]$, $i = 1, 2, \ldots$ are disjoint, it follows that $i \delta \rho\left(\frac{\varepsilon}{2}\right) \leq \int_0^{t_i} \rho(f(s)) ds$ for all $i = 1, 2, \ldots$, which contradicts the fact that $\int_0^{+\infty} \rho(f(t)) dt < +\infty$.

Now consider the case when $f$ is bounded and let $R := \sup\{f(t) : t \geq 0\}$ and $c := \min\left\{\rho(s) : \frac{\varepsilon}{2} \leq s \leq R\right\}$. It follows from (5.15) that $\rho(f(s)) \geq c$ for all $s \in [t_i - \delta, t_i]$, $i = 1, 2, \ldots$. Since the intervals $[t_i - \delta, t_i]$, $i = 1, 2, \ldots$ are disjoint, it follows that $i \delta c \leq \int_0^{t_i} \rho(f(s)) ds$ for all $i = 1, 2, \ldots$, which contradicts the fact that $\int_0^{+\infty} \rho(f(t)) dt < +\infty$.

Case 2: $-f$ is quasi-uniformly continuous

Since $-f$ is a quasi-uniformly continuous function there exists $\delta > 0$ such that $f(t) - f(\tau) > -\varepsilon/2$ for all $t \in [\tau, \tau + \delta]$.

Without loss of generality (by taking a subsequence if necessary) we may assume that $t_1 \geq 1 + \delta$ and $t_{i+1} \geq t_i + \delta + 1$ for all $i = 1, 2, \ldots$. Notice that

$$\frac{\varepsilon}{2} \leq f(t_i) - \frac{\varepsilon}{2} < f(t) \text{ for all } t \in [t_i, t_i + \delta], \ i = 1, 2, \ldots. \tag{5.16}$$

If $\rho : \Re_+ \to \Re_+$ is non-decreasing then (5.16) implies $\rho\left(\frac{\varepsilon}{2}\right) \leq \rho(f(t))$ for all $t \in [t_i, t_i + \delta]$, $i = 1, 2, \ldots$. Since the intervals $[t_i, t_i + \delta]$, $i = 1, 2, \ldots$ are disjoint, it follows that $(i-1) \delta \rho\left(\frac{\varepsilon}{2}\right) \leq \int_0^{t_i} \rho(f(s)) ds$ for all $i = 1, 2, \ldots$, which contradicts the fact that $\int_0^{+\infty} \rho(f(t)) dt < +\infty$.

Now consider the case when $f$ is bounded and let $R := \sup\{f(t) : t \geq 0\}$ and $c := \min\left\{\rho(s) : \frac{\varepsilon}{2} \leq s \leq R\right\}$. It follows from (5.16) that $\rho(f(t)) \geq c$ for all $t \in [t_i, t_i + \delta]$, $i = 1, 2, \ldots$. Since the intervals $[t_i, t_i + \delta]$, $i = 1, 2, \ldots$ are disjoint, it follows that $(i-1) \delta c \leq \int_0^{t_i} \rho(f(s)) ds$ for all $i = 1, 2, \ldots$, which contradicts the fact that $\int_0^{+\infty} \rho(f(t)) dt < +\infty$. ◁

## 5.6. Proof of Theorem 2

Proposition 1 guarantees Lagrange output stability on $\Omega$ and Lyapunov output stability on $\Omega$.

Notice that (2.6) implies differential inequality (5.5) for all $x_0 \in \Omega$. Consequently, (5.5) implies that the following estimate holds for all $x_0 \in \Omega$ and $t \geq 0$:

$$V(\phi(t, x_0)) \leq V(x_0) \tag{5.17}$$

Inequalities (2.9), (2.10) imply that one of the differential inequalities holds for all $x_0 \in \Omega$:

$$\frac{d}{dt} W(\phi(t, x_0)) \leq \gamma(V(\phi(t, x_0))), \text{ for all } t \geq 0 \tag{5.18}$$

or

$$\frac{d}{dt} W(\phi(t, x_0)) \geq -\gamma(V(\phi(t, x_0))), \text{ for all } t \geq 0 \tag{5.19}$$

Taking into account inequality (5.17) and defining

$$G(x_0) = \max\{\gamma(s) : 0 \leq s \leq V(x_0)\}, \text{ for all } x_0 \in \Omega \tag{5.20}$$

we obtain from (5.18) and (5.19) that one of the following differential inequalities holds for all $x_0 \in \Omega$:

$$\frac{d}{dt} W(\phi(t, x_0)) \leq G(x_0), \text{ for all } t \geq 0 \tag{5.21}$$

or

$$\frac{d}{dt} W(\phi(t, x_0)) \geq -G(x_0), \text{ for all } t \geq 0 \tag{5.22}$$

Proposition 2 in conjunction with (5.21) and (5.22) implies that either for all $x_0 \in \Omega$ the mapping $t \to W(\phi(t, x_0))$ is quasi-uniformly continuous (when (5.21) holds) or for all $x_0 \in \Omega$ the mapping $t \to -W(\phi(t, x_0))$ is quasi-uniformly continuous (when (5.22) holds).

Using (5.5) we obtain (5.12), which shows that for every $x_0 \in \Omega$ it holds that $\int_0^{+\infty} \rho(W(\phi(s, x_0))) ds \leq V(x_0)$.

If $\rho$ is non-decreasing then the application of Lemma 3 to the function $f(t) = W(\phi(t, x_0))$ shows that $\lim_{t \to +\infty}(W(\phi(t, x_0))) = 0$ for every $x_0 \in \Omega$. The fact that $\lim_{t \to +\infty}(y(t, x)) = 0$ is a direct consequence of inequality (2.7) and the fact that $\lim_{t \to +\infty}(W(\phi(t, x))) = 0$.

On the other hand, if (2.11) holds then we get for all $x_0 \in \Omega$:

$$W(\phi(t, x_0)) \leq \zeta(V(\phi(t, x_0))), \text{ for all } t \geq 0 \tag{5.23}$$

Taking into account inequality (5.17) and defining

$$Z(x_0) = \max\{\zeta(s): 0 \leq s \leq V(x_0)\}, \text{ for all } x_0 \in \Omega \tag{5.24}$$

we obtain from (5.23) that the following estimate holds for all $x_0 \in \Omega$:

$$W(\phi(t, x_0)) \leq Z(x_0), \text{ for all } t \geq 0 \tag{5.25}$$

Estimate (5.25) shows that for all $x_0 \in \Omega$ the mapping $t \to W(\phi(t, x_0))$ is bounded. Therefore, the application of Lemma 3 to the function $f(t) = W(\phi(t, x_0))$ shows that $\lim_{t \to +\infty}(W(\phi(t, x_0))) = 0$ for every $x_0 \in \Omega$. The fact that $\lim_{t \to +\infty}(y(t, x)) = 0$ is a direct consequence of inequality (2.7) and the fact that $\lim_{t \to +\infty}(W(\phi(t, x))) = 0$. ◁

5.7. Proof of Theorem 3

Using (3.2), (3.9) and the definition $V(x) := P(y) + \frac{\gamma}{2}|z|^2$ we can establish the inequality

$$\frac{d}{dt}V(x) \leq -Q(y), \text{ for all } x = (y, z) \in \Re^n \times \Re^p \tag{5.26}$$

Since $Q$ and $P$ are positive definite functions and since $P$ is radially unbounded, there exists a continuous, positive definite function $\rho: \Re_+ \to \Re_+$ with the following property (see Proposition 2.2 on page 107 in [6]):

$$Q(y) \geq \rho(P(y)), \text{ for all } y \in \Re^n \tag{5.27}$$

Moreover, $V(x) := P(y) + \frac{\gamma}{2}|z|^2$ is radially unbounded. Therefore, there exist functions $a, b \in K_\infty$ with the following properties (see Proposition 2.2 on page 107 in [6]):

$$a(|y|) \leq P(y), \text{ for all } y \in \Re^n \tag{5.28}$$

$$V(x) = P(y) + \frac{\gamma}{2}|z|^2 \leq b(|x|), \text{ for all } x = (y, z) \in \Re^n \times \Re^p \tag{5.29}$$

Let $x_0 \in \Re^n \times \Re^p$ be given. The initial-value problem (3.9) with initial condition $x(0) = x_0$ has a unique solution $x:[0, t_{\max}) \to \Re^n \times \Re^p$, where $t_{\max} \in (0, +\infty]$ is the maximal existence time of the corresponding solution. The differential inequality (5.26) implies that

$$V(x(t)) \leq V(x_0), \text{ for all } t \in [0, t_{\max}) \tag{5.30}$$

Define the set $\Psi := \{x \in \Re^n \times \Re^p : V(x) \leq V(x_0)\}$. Since $V(x) := P(y) + \frac{\gamma}{2}|z|^2$ is a positive definite and radially unbounded function, $\Psi$ is a compact set. Since estimate (5.30) implies that $x(t) \in \Psi$ for all $t \in [0, t_{\max})$, it follows that $t_{\max}$ cannot be finite. Consequently, $t_{\max} = +\infty$. Estimate (5.30) shows that the set $\Psi := \{ x \in \Re^n \times \Re^p : V(x) \leq V(x_0) \}$ is a positively invariant set for the dynamical system (3.9). Similarly, the set $\Omega$ defined by (3.10) is a positively invariant set for the dynamical system (3.9).

Theorem 2 (and Remark 2.2(ii)) in conjunction with the fact that $V(x) := P(y) + \frac{\gamma}{2}|z|^2$ is a radially unbounded function implies that system (3.9) with output (3.5) is AOS on $\Re^n \times \Re^p$.

If $x = (y, z) \in \Omega$ then definition (3.10) implies that $|z| \leq \sqrt{2L}$. Consequently, we obtain for all $x = (y, z) \in \Omega$:

$$\left| \nabla P(y) g(y) (\varphi(y))' z \right| \leq \left| \nabla P(y) g(y) \right| \left| \varphi(y) \right| |z|$$
$$\leq \sqrt{2L} \left| \nabla P(y) g(y) \right| \left| \varphi(y) \right|$$
$$\leq L \mu(y) \left| \nabla P(y) g(y) \right|^2 + \frac{1}{2\mu(y)} \left| \varphi(y) \right|^2$$

It follows from (3.7) and the above inequality that the following inequality holds for all $x = (y, z) \in \Omega$:

$$\left| \nabla P(y) g(y) (\varphi(y))' z \right| \leq L \mu(y) \left| \nabla P(y) g(y) \right|^2 + \frac{1}{2} Q(y) \tag{5.31}$$

Using inequalities (3.2) and (5.31), we get from (3.9) for all $x = (y, z) \in \Omega$:

$$\frac{d}{dt} P(y) = \nabla P(y) \dot{y}$$
$$= \nabla P(y) f(y) + \nabla P(y) g(y) k(y) - L \mu(y) \left| \nabla P(y) g(y) \right|^2 - \nabla P(y) g(y) (\varphi(y))' z \tag{5.32}$$
$$\leq -Q(y) - L \mu(y) \left| \nabla P(y) g(y) \right|^2 + \left| \nabla P(y) g(y) (\varphi(y))' z \right| \leq -\frac{1}{2} Q(y) \leq 0$$

Inequalities (5.26), (5.27), (5.28), (5.32) show that all assumptions of Theorem 1 hold with $W(x) = P(y)$ and UAOS on $W$ follows. ◁

## 6. Concluding Remarks

We have provided a Lyapunov-based condition to ensure uniform output asymptotic stability in the case when the Lyapunov function does not involve the Lyapunov function itself, but rather a function $W$ which can only vanish when the output is zero. Unlike classical tools such as Barbălat's lemma and Barbashin-Krasovskii-LaSalle's invariance principle, the convergence of the output to zero is guaranteed to be uniform on bounded sets of the state provided that this function $W$ does not increase along the system's solutions. We have shown that this result can be useful in the analysis of control systems under adaptive control.

We have also presented a relaxation of Barbălat's lemma, which does not require uniform continuity of the considered function. This requirement is replaced by "quasi-uniform continuity", which can be inferred by showing that the derivative of the considered function is upper bounded (but not necessarily lower-bounded). We have shown through academic examples that this relaxation can prove useful in situations where Barbălat's lemma cannot be applied.

We have also shown that these results are not confined to finite-dimensional systems, but also hold in infinite dimension, thus covering time-delay systems and systems ruled by partial differential equations.

One of the advantages of Barbălat's lemma is its ability to be applied to time-varying systems: future work should aim at extending the results presented here to such class of systems, which could have practical interest for instance in adaptive control for trajectory tracking. Further investigations are also needed to cover control systems with uncertain parameters that do not satisfy the matching condition. Finally, uniform stability properties own the advantage of providing robustness to exogenous disturbances or model imprecision: future work is needed to formalize this fact in a general context.